\documentclass{article}
\usepackage{amssymb,amsfonts}

\newcommand{\trace}{\mathop{\rm Tr}\nolimits}

\newcommand{\twomat}[4]{\left(\begin{array}{cc}#1&#2\\#3&#4\end{array}\right)}

\newcommand{\schatten}[2]{\left|\left|\,{#2}\,\right|\right|_{#1}}

\DeclareRobustCommand\openone{\leavevmode\hbox{\small1\normalsize\kern-.33em1}}
\newcommand{\id}{\mathrm{\openone}}
\newcommand{\pf}{\textit{Proof.}}
\newcommand{\eop}{\hfill$\square$\par\vskip24pt}
\newcommand{\be}{\begin{equation}}
\newcommand{\ee}{\end{equation}}
\newcommand{\bea}{\begin{eqnarray}}
\newcommand{\eea}{\end{eqnarray}}
\newcommand{\beas}{\begin{eqnarray*}}
\newcommand{\eeas}{\end{eqnarray*}}

\newtheorem{theorem}{Theorem}
\newtheorem{lemma}{Lemma}

\newtheorem{proposition}{Proposition}
\newcount\minute
\newcount\hour
\def\currenttime{%
    \minute\time
    \hour\minute
    \divide\hour60
    \the\hour:\multiply\hour60\advance\minute-\hour\the\minute}
\begin{document}
\title{On the Araki-Lieb-Thirring inequality}
\author{Koenraad M.R. Audenaert \\
Institute for Mathematical Sciences, Imperial College London \\
53 Prince's Gate, London SW7 2PG, United Kingdom}
\date{\today, \currenttime}
\maketitle

\begin{abstract}
We prove an inequality that complements the famous Araki-Lieb-Thirring (ALT) inequality
for positive matrices $A$ and $B$, by giving a
lower bound on the quantity $\trace[A^r B^r A^r]^q$ in terms of
$\trace[ABA]^{rq}$ for $0\le r\le 1$ and $q\ge0$, whereas the ALT inequality gives an upper bound.
The bound contains certain norms of $A$ and $B$ as additional ingredients and is therefore of a different nature than the
Kantorovich type inequality obtained by Bourin (\textit{Math.\ Inequal.\ Appl.} \textbf{8}(2005) pp.\ 373--378) and others.
Secondly, we also prove a generalisation of the ALT inequality to general matrices.
\end{abstract}

\section{Introduction\label{sec1}}
A famous inequality with a lot of applications in mathematics and mathematical physics
is the Araki-Lieb-Thirring inequality \cite{araki,lieb}:
\begin{theorem}[Araki-Lieb-Thirring]
For $A,B\ge0$, $q\ge0$, and for $0\le r\le 1$, the following inequality holds:
\be\label{eq:ALT}
\trace[A^r B^r A^r]^q \le \trace[ABA]^{rq},
\ee
while for $r\ge1$, the inequality is reversed.
\end{theorem}
Together with its companion, the Lieb-Thirring inequality, it has been extended in various directions, see for example \cite{wang}
and references therein.

In this paper we do two things. In Section \ref{sec:reverse} we obtain complementary inequalities. That is, for $0\le r\le1$ we obtain
upper bounds on $\trace[ABA]^{rq}$ (in terms of the quantity $\trace[A^r B^r A^r]^q$), and lower bounds for $r\ge1$.
These bounds contain certain norms of $A$ and $B$ as additional ingredients, and are therefore of a different nature than the
Kantorovich type inequalities obtained by Bourin \cite{bourin} and others.
Second, in Section \ref{sec:gen}, we find a generalisation of the ALT inequality to general matrices.
\section{A complementary inequality\label{sec:reverse}}
In this Section we want to obtain upper bounds on $\trace[ABA]^{rq}$ in terms of the quantity $\trace[A^r B^r A^r]^q$
(for $0\le r\le1$).
A remark one can make right away
is that both quantities $\trace[ABA]^{rq}$ and $\trace[A^r B^r A^r]^q$
have the same degrees of homogeneity in $A$ and $B$. Any upper bound on $\trace[ABA]^{rq}$ that only depends on $\trace[A^r B^r A^r]^q$
should therefore be linear in the latter quantity. Unfortunately, numerical calculations show that the required proportionality factor
should be infinitely large to accomodate all possible $A$ and $B$. This means that extra ingredients are needed to obtain
a reasonable upper bound.

One way to do this is to supply the values of the extremal eigenvalues of $B$,
yielding Kantorovich-type inequalities.
This has been investigated by Bourin in \cite{bourin}, who obtained
the inequalities
$$
K(a,b,r)^{-1}\lambda^\downarrow((ABA)^r)\le \lambda^\downarrow(A^r B^r A^r)
\le K(a,b,r)\lambda^\downarrow((ABA)^r),
$$
for $r\ge1$, $A\ge0$ and $0<b\le B\le a$,
where $K(a,b,r)$ is the Ky Fan constant
$$
K(a,b,r):=\frac{a^r b-a b^r}{(r-1)(a-b)}\left(\frac{r-1}{r}\,\,\frac{a^r-b^r}{a^r b-a b^r}\right)^r,
$$
and where $\lambda^\downarrow(A)$ denotes the vector of eigenvalues of $A$ sorted in non-increasing order.
Previous results in this direction were obtained by
Furuta \cite{furuta} and Fujii, Seo and Tominaga \cite{fujii}.

In this work, we have followed a different route
and have found an upper bound by including norms of both $A$ and $B$ as additional ingredients.
In fact, we have found a whole family of such bounds.
The simplest bound in this family, but also the weakest, is given by
\begin{proposition}
Let $A,B\ge0$. For $q\ge0$ and $r\ge0$, we have the upper bound
\be\label{eq:water}
\trace[ABA]^{rq} \le ||A||^{2rq} \trace B^{rq},
\ee
while for $q\le0$ the inequality is reversed.
\end{proposition}
Here, $||.||$ denotes the operator norm, which for positive semidefinite (PSD) matrices is nothing but the
largest eigenvalue.

\pf
Put $p=rq\ge0$.
Note first that
$$
\trace[ABA]^p = \trace[B^{1/2}A^2B^{1/2}]^p.
$$
From the basic inequality $A^2\le ||A^2||\id = ||A||^2\id$ follows
$$
B^{1/2}A^2B^{1/2} \le ||A||^2 B.
$$
If $p$ is between 0 and 1, we may take the $p$-th power of both sides ($x\mapsto x^p$ is then operator monotone).
Taking the trace of both sides then yields (\ref{eq:water}).
If $p$ is larger than 1, we may take the Schatten $p$-norm of both sides, by Weyl-monotonicity of unitarily
invariant (UI) norms.
Taking the $p$-th power of both sides again yields (\ref{eq:water}).
\eop

For reasons that will immediately become clear,
we call this inequality the ``water''-inequality, to express the fact that it is rather weak, and not very spiritual.
In contrast, we call the Araki-Lieb-Thirring inequality the ``wine''-inequality,
because it is too strong for our purposes:
it gives a lower bound, rather than an upper bound.

We can obtain better upper bounds by ``cutting the wine with the water''. Fixing $A$ and $B$,
some $t$ obviously must exist, with $0\le t\le 1$, such that the following holds:
\beas
\trace[ABA]^{rq} &\le& (\mbox{water})^t \,\,(\mbox{wine})^{1-t} \\
&=& \left(||A||^{2rq} \trace B^{rq}\right)^t\,\,\left(\trace[A^r B^r A^r]^q\right)^{1-t}.
\eeas
Of course, this would be a rather pointless (and disappointing)
exercise if there were some $A$ and $B$ for which the smallest
valid value of $t$ would be 1, because then inequality (\ref{eq:water}) would be the only
upper bound valid for all $A$ and $B$.
Fortunately, numerical experiments revealed the fact (which we will prove below) that here
any value of $t$ between $1-r$ and $1$ yields an upper bound, for \textit{any} $A$ and $B$.
This yields the promised family of inequalities, of which the sharpest and most
relevant one is the one with $t=1-r$.
%
\begin{theorem}\label{th:reverse}
For $A,B\ge0$, $q\ge0$ and $0\le r\le 1$,
\be\label{eq:waterwine}
\trace[ABA]^{rq} \le \left(||A||^{2rq} \trace B^{rq}\right)^{1-r}\,\,\left(\trace[A^r B^r A^r]^q\right)^r.
\ee
For $r\ge 1$, the inequality is reversed.
\end{theorem}
This inequality is sharp, just like the original ALT inequality, as can be seen by taking scalar $A$ and $B$.

An equivalent formulation of inequality (\ref{eq:waterwine}) is
\be\label{eq:waterwine2}
||(ABA)^r||_q \le \left(||A||^{2r} ||B^r||_q\right)^{1-r}\,\,||A^r B^r A^r||_q^r.
\ee
Yet another formulation is obtained if one notes the equality
$$
||(ABA)^r||_q = ||ABA||^r_{rq} = ||B^{1/2}A^2 B^{1/2}||^r_{rq} = ||(B^{1/2}A^2 B^{1/2})^r||_q,
$$
by which we get
\bea
||(ABA)^r||_q &=& ||(B^{1/2}A^2 B^{1/2})^r||_q \nonumber\\
&\le& \left(||B^{1/2}||^{2r} ||A^{2r}||_q\right)^{1-r}\,\,||B^{r/2} A^{2r} B^{r/2}||_q^r \nonumber\\
&=& \left(||B||^{r} ||A^{2r}||_{q}\right)^{1-r}\,\,||A^r B^r A^r||_q^r. \label{eq:waterwine3}
\eea

It is this last formulation that we will consider in the following proof.

\pf
Let first $0\le r\le 1$.
Then, for $X\ge0$, $||X^{1-r}|| = ||X||^{1-r}$. 
Since $B\ge0$ we can write $B=B^r B^{1-r} \le ||B||^{1-r} B^r$.
Thus
\be
ABA \le ||B||^{1-r} A B^r A = ||B||^{1-r} A^{1-r} (A^r B^r A^r) A^{1-r}.\label{eq:aba}
\ee

Consider first the easiest case $q=\infty$.
Taking the operator norm of both sides of (\ref{eq:aba}) then gives
$$
||ABA|| \le ||B||^{1-r} ||A^{1-r} (A^r B^r A^r) A^{1-r}||
\le ||B||^{1-r} ||A^{1-r}||^2 ||A^r B^r A^r||,
$$
where the last inequality follows from submultiplicativity of the operator norm.
One obtains (\ref{eq:waterwine3}) for $q=\infty$
by taking the $r$-th power of both sides, and noting (again) $||X||^r = ||X^r||$.

\bigskip

To prove (\ref{eq:waterwine3}) for general $q$, let us first take the $r$-th power of both sides of (\ref{eq:aba}),
which preserves the ordering because of operator monotonicity of $x\mapsto x^r$ for $0\le r\le1$:
$$
(ABA)^r \le ||B^r||^{1-r} \left(A^{1-r} (A^r B^r A^r) A^{1-r}\right)^r.
$$
Thus, on taking the $q$-norm of both sides (or $q$-\textit{quasi}norm if $0<q<1$),
\be
||(ABA)^r||_q
\le ||B^r||^{1-r} ||\left(A^{1-r} (A^r B^r A^r) A^{1-r}\right)^r||_q. \label{eq:abar}
\ee

We can now apply a generalisation of H\"older's inequality, by which
for all positive real numbers $s,t,u$ such that $1/s+1/t=1/u$ we have (\cite{bhatia}, Eq.\ (IV.43))
$$
|||\,\,|XY|^u\,\,|||^{1/u} \le |||\,\,|X|^s\,\,|||^{1/s} \,\, |||\,\,|Y|^t\,\,|||^{1/t},
$$
for all $X,Y$ and for all UI norms $|||.|||$.
In fact, this inequality extends to UI quasinorms like the Schatten $q$-quasinorms for $0<q<1$.

For $X,Y\ge0$, two successive applications of this inequality yield
$$
|||\,(XYX)^u\,|||^{1/u} \le |||X^{2s}|||^{1/s} \,\, |||Y^t|||^{1/t}.
$$
We apply the latter inequality to the second factor of the right-hand side (RHS) of (\ref{eq:abar}), with
the substitutions $X=A^{1-r}$, $Y=A^r B^r A^r$, $u=r$, $s=r/(1-r)$ (the positivity of which requires $r$ to lie between 0 and 1),
$t=1$, and $|||.||| = ||.||_q$,
giving
$$
||\left(A^{1-r} (A^r B^r A^r) A^{1-r}\right)^r||_q^{1/r}
\le ||A^{2r}||_q^{(1-r)/r} \,\,||A^r B^r A^r||_q.
$$
Taking the $r$-th power of both sides and substituting in (\ref{eq:abar})
yields (\ref{eq:waterwine3}).

\bigskip

The case $r\ge1$ follows very easily from the case $0\le r\le 1$ by making in (\ref{eq:waterwine}) the substitutions
$A'=A^r$, $B'=B^r$, $r'=1/r$, $q'=qr$. Taking the $r'$-th power of both sides, rearranging factors, and subsequently dropping primes
yields (\ref{eq:waterwine}) for $r\ge1$.
\eop

As a special case of Theorem \ref{th:reverse} we consider the comparison between $\trace[AB]$ and $||AB||_1$, for $A,B\ge0$.
The lower bound $||AB||_1\ge \trace[AB]$ is well-known and easy to prove.
Indeed, since $A^{1/2}BA^{1/2}$ is normal,
\beas
\trace[AB]&=&\trace[A^{1/2}BA^{1/2}] \\
&=& ||A^{1/2}BA^{1/2}||_1 \\
&\le& ||A^{1/2}A^{1/2}B||_1=||AB||_1.
\eeas
To obtain an upper bound, (\ref{eq:waterwine}) with $r=1/2$ and $q=1$ yields
\bea
||AB||_1 &=& \trace(AB^2A)^{1/2} \nonumber \\
&\le& \left(||A||\trace B\,\,\trace(A^{1/2}BA^{1/2})\right)^{1/2} \nonumber \\
&=& \left(||A||\trace B\,\,\trace(AB)\right)^{1/2}.
\eea
\section{Generalisations of the Araki-Lieb-Thirring inequality\label{sec:gen}}
In this Section we want to find a generalisation of the ALT inequality to general matrices
$A$ and $B$.
As a first step, we generalise it to the case of general $A$ while keeping $B\ge0$.
\begin{proposition}\label{prop:ALT2}
For any matrix $A$ and $B\ge0$, for $q\ge1$, and for any UI norm,
\be
|||(ABA^*)^q|||\le |||\,\, |A|^q B^q |A|^q\,\, |||.
\ee
\end{proposition}

\pf
Let the polar decomposition of $A$ be $A=U|A|$, where $|A|=(A^*A)^{1/2}$ is the modulus of $A$ and $U$ is unitary.
Then, for any UI norm,
$$
|||(ABA^*)^q|||
= ||| (U|A|\,\, B\,\, |A|U^*)^q |||
= ||| (|A|\,\, B\,\, |A|)^q |||
\le |||\,\, |A|^q B^q |A|^q \,\,|||,
$$
where in the last step we used the ALT inequality proper.
\eop

The second step is to generalise this statement to the case where $B$ is Hermitian.
First we need a Lemma.
\begin{lemma}
For $X,Y\ge0$, and any UI norm, $|||X-Y|||\le |||X+Y|||$.
\end{lemma}

\pf
The matrix $\twomat{X}{0}{0}{Y}$ is unitarily equivalent with
$\twomat{X+Y}{X-Y}{X-Y}{X+Y}/2$. Since $X,Y\ge0$, the latter matrix is also PSD, whence there exists a contraction $K$
such that $X-Y=(X+Y)^{1/2}K(X+Y)^{1/2}$. Moreover, as $X-Y$ and $X+Y$ are Hermitian, so is $K$.
Therefore
$$
|||X-Y||| = |||(X+Y)^{1/2}K(X+Y)^{1/2}||| \le |||(X+Y)K|||,
$$
where we have used the fact that $|||AB|||\le |||BA|||$ whenever $AB$ is normal.
Every contraction $K$ can be written as a convex combination of unitaries, so by convexity
of norms, and the fact that we're considering UI norms, we have $|||(X+Y)K|||\le |||X+Y|||$.
\eop

\begin{proposition}\label{prop:LTH}
For any matrix $A$ and Hermitian $B$, for $q\ge 1$, and for any UI norm,
\be
|||\,\,|ABA^*|^q\,\,||| \le |||\,\,|A|^q |B|^q |A|^q \,\,|||.
\ee
\end{proposition}

\pf
Let the Jordan decomposition of $B$ be $B=B^+ - B^-$, where $B^+$ and $B^-$ (the positive and negative part, respectively) are both PSD.
Then the two terms in the right-hand side of $ABA^* = AB^+A^* - AB^-A^*$ are also PSD.
Thus
$$
|||\,\,|ABA^*|^q\,\,||| \le |||(AB^+A^* + AB^-A^*)^q|||.
$$
This follows from the Lemma applied to the norm $|||\,\,|\cdot|^q\,\, |||^{1/q}$.

But as $AB^+A^* + AB^-A^* = A|B|A^*$, we get
$$
|||\,\,|ABA^*|^q\,\,||| \le |||(A\,|B|\,A^*)^q|||\le |||\,\,|A|^q |B|^q |A|^q \,\,|||,
$$
where in the last step we used Proposition \ref{prop:ALT2}.
\eop

If we now specialise to Schatten $p$-norms,
we can drop the conditions on $B$:
\begin{theorem}\label{prop:LTG}
For general matrices $A$ and $B$, and for $p,q\ge 1$,
\be
\schatten{p}{|ABA^*|^q} \le \schatten{p}{|A|^q \frac{|B|^q+|B^*|^q}{2}|A|^q}.
\ee
\end{theorem}

\pf
We use Proposition \ref{prop:LTH} with
$A=\twomat{A'}{0}{0}{A'}$ and $B=\twomat{0}{B'}{B'^*}{0}$.
Noting that
$$
\left|\twomat{0}{X}{X^*}{0}\right| = \twomat{|X^*|}{0}{0}{|X|},
$$
this yields, after dropping primes,
$$
\schatten{p}{\twomat{|AB^*A^*|^q}{0}{0}{|ABA^*|^q}} \le
\schatten{p}{\twomat{|A|^q|B^*|^q|A|^q}{0}{0}{|A|^q |B|^q |A|^q)}}.
$$
Using the facts $||X\oplus Y||_p^p = ||X||_p^p + ||Y||_p^p$
and $||\,\,|X|^q\,\,||_p = ||\,\,|X^*|^q\,\,||_p$,
yields the statement of the Theorem.
\eop

\section*{Acknowledgments}
This work was supported by the Institute for Mathematical Sciences, Imperial College London,
and is part of the QIP-IRC (www.qipirc.org) supported by EPSRC (GR/S82176/0).
I thank an anonymous referee for pointing out reference \cite{bourin}.

\end{document}